*Nonlinear retracts and the geometry of Banach spaces*

M. A. Sofi

*Abstract:* In the nonlinear geometry of Banach spaces where the objects in the category are Banach spaces as in the linear case, the morphisms in the new setting are taken to comprise of certain nonlinear maps involving say, Lipschitz maps and, in some cases, uniformly continuous or coarse (resp. coarse Lipschitz) mappings arising from the underlying metric (resp. uniform) structure attached to the norm of the Banach space. The question as to what extent the Lipschitz (resp. uniform) structure may be used to capture the full linear structure of a Banach space has been one of the most fundamental problems pertaining to the nonlinear structure of Banach spaces since this line of investigation was undertaken by Lindenstrauss in the early seventies. This line of research which is subsumed under the so-called Ribe program broadly underscores the view that metric spaces encode a much deeper and hidden structure than is apparent. It's truly surprising how the (linear) structure of a Banach space gets captured to a considerable extent by its metric space structure. This point of view has led to deep insights into Banach space theory that has paved the way for these ideas being employed in seemingly unrelated disciplines including harmonic analysis, geometric group theory, Riemannian geometry, theoretical computer science, probability theory, data structures and many other domains of mathematics.

*1. Background*

The notion of retracts is a generic theme in topology and analysis. In the theory of Banach spaces, the idea of a Lipschitz/uniform retract may be seen as a nonlinear counterpart of a complemented subspace of a Banach space which has been widely studied in the linear theory of Banach spaces. In what follows we shall see how the structure of Lipschitz retracts captures the linear structure of Banach spaces, which in most cases mimics the theory of complemented subspaces of Banach spaces. An important example of this phenomenon, to be proved in Section 3 is provided by the Lipschitz analogue of the celebrated Lindenstrauss-Tzafiriri complemented subspaces problem characterising isomorphs of Hilbert spaces. For an advanced introduction to Banach space theory, we refer to [1] whereas for a comprehensive treatment of Lipschitz mappings we refer to the encyclopaedic treatise [10] and of their role in the nonlinear geometry of Banach spaces, the best source of information is the authoritative monograph[4]. Another useful source is [45] providing a thorough description of the theory underlying the space $Lip_0(M)$ as an algebra and its canonical predual, the free space $\mathfrak{I}(X)$ of X.

The main aim of this note is to give an introduction and overview of recent developments surrounding the nonlinear geometry of Banach spaces including the present state of art of the theory. Barring a few instances where a brief outline of the proofs has been indicated, most of the results which appeal to a much deeper and broader understanding of the technical machinery underlying the theory have been presented without proofs.

As we shall soon learn (see Section 3 and 4), an interesting feature of Lipschitz geometry of Banach spaces is the crucial role that the separability of Banach spaces plays in our considerations involving the proofs of some deep theorems. On the other side, it turns out that almost all the highly nontrivial counterexamples to certain natural questions arising in the Lipschitz classification of Banach spaces have been constructed for the nonseparable case. The construction of such counterexamples in the separable case remains a pesky issue as the main obstacle in these constructions is the huge bulk of important positive results in the theory which are valid under the assumption of separability of Banach spaces in question. An important instance of such results is provided by the well-known Godefroy-Kalton theorem on the existence of a right Lipschitz inverse for a given bounded linear mapping between Banach spaces (See Section 5)

Throughout, we shall use X, Y, Z to denote Banach spaces, unless otherwise stated whereas $X^*$ shall denote the dual of X. The closed unit ball of X shall be denoted by $B_X$: $B_X = \{x \in X, \|x\| \leq 1\}$. The symbol L(X, Y) shall be used for the Banach space of bounded linear operators from X into Y. The following notation shall also be used in the course of our discussion:

$$\ell_p^n = (\mathbb{R}^n, \|\ \|_p), where \|(x_i)\|_p = \left(\sum_{i=1}^n |x_i|^p\right)^{1/p}, 1 \leq p < \infty$$

$$\ell_\infty^n = (\mathbb{R}^n, \|\ \|_\infty). \text{ Here, } (\|(x_i)\|_\infty = \max_{1 \leq i \leq n} |x_i|)$$

$$\ell_p = \left\{(x_i); \|(x_i)\|_p = \left(\sum_{i=1}^\infty |x_i|^p\right)^{1/p} < \infty\right\}, 1 \leq p < \infty.$$

$\ell_\infty$, Banach space of all bounded sequences equipped with the norm $(\|(x_i)\|_\infty = \sup_{i \geq 1}|(x_i)|)$.

$c_0$, the (Closed) subspace of $\ell_\infty$ consisting of null sequences (with the induced norm).

Definition 1.1: Given a set A, $B \subset A$ and a mapping $f: A \to B$, f is said to be a *retraction* on B if $f(x) = x, \forall\ x \in B$. In this case, B is called a *retract* of A. Continuous, uniform and Lipschitz retracts may be defined analogously, depending upon whether f is a continuous, uniformly continuous or a Lipschitz mapping acting between topological (resp. metric) spaces defined below.

Definition 1.2: A mapping $f: M \to N$ acting between metric spaces M and N is said to be a Lipschitz map if the following holds:

$$d(f(x), f(y)) \leq k d(x,y), x, y \in M, k \geq 1.$$

For k=1, we say that $f$ is a *nonexpansive* (NE) mapping. We begin with some examples in the finite dimensional setting.

Theorem 1.3 (L. E. J. Brouwer): $S^n$ is not a (continuous) retract of $B^n$. Equivalently, every continuous self-map on $B^n$ admits a fixed point.

Example 1.4: The closed unit disc (more generally every closed and convex subset) of the Euclidean plane is a nonexpansive retract (NR).

More generally, we have

Theorem 1.5: Every closed convex subset of the Minkowski plane is an (NR).

In higher dimensions, we have

Example 1.6: Not every closed, convex subset of the 3-dimensional space is an (NR). Consider the plane P: x+y+z=4 in $\ell_\infty^3$. Then the point D(1,1,1), which is distant 1 from *each* of the points A(2,2,0), B(2,0,2) and C(0,2,2) on P shall get mapped to some point on P which would necessarily bear a distance strictly larger than 1 *at least from one* of the points A, B and C. This also shows that the Lipschitz constant of an extension of a Lipschitz map may not be preserved.

Proposition 1.7: Closed convex subsets of $\ell_2^n$ are (NR).
Corollary 1.8: For X with $\dim X < \infty$, every closed convex subset of X is a Lipschitz retract.

## 2. Infinite dimensional setting

The following example (attributed to S. Kakutani) shows that the infinite dimensional analogue of Brouwer's fixed point theorem does not hold.

Ex 2.1: The continuous map $f: B_{\ell_2} \to B_{\ell_2}$ given by:

$$f(x) = \left(\sqrt{1 - \|x\|^2}, x_1, x_2, , , \ldots x_n, \quad \ldots \right)$$

is a retraction of $B_{\ell_2}$ onto the unit sphere.

The question regarding the validity of the above phenomenon in general Banach spaces is treated in the following theorem of W. A. Kirk ([4], Theorem 3.3. See also [15], Chapter 4) which shows that the above phenomenon holds in all infinite dimensional Banach spaces.

Theorem 2.2: Every infinite dimensional Banach space admits a fixed point free continuous map on its unit ball. Equivalently, the unit sphere of an infinite dimensional Banach space is a (continuous) retract of the unit ball. More generally, every noncompact closed convex subset of an infinite dimensional Banach space X is a retract of X.

The following result is closely related to the previous theorem.

Proposition 2.3: Every closed convex subset A of a Banach space is an absolute (continuous) retract, i.e., A is a (continuous) retract of every metric space containing it as a closed subset.

The above result is a consequence of Michael's selection theorem given below:

Theorem 2.4 ([4], Chapter 1): Let H(X) denote the collection of all nonempty, closed, convex subsets of a Banach space X equipped with the Hausdorff metric:

$$d_H(A, B) = sup(\{d(x, A); x \in B\} \cup \{d(x, B); x \in A\}).$$

Further, let $\phi: M \to H(X)$ be a (set-valued) mapping which is *lower semi-continuous* (l.s.c.): $\forall \mathcal{U}$, an open subset of X, the set $\{x \in M; \phi(x) \cap \mathcal{U} \neq \emptyset\}$ is an open subset of M. Then $\phi$ admits a continuous selection $\varphi: M \to X$; $\varphi(x) \in \phi(x), for\ every\ x \in M$. In fact, it is possible to show that there exists a continuous *selector* $\psi: H(X) \to X$, i.e., $\psi(A) \in A, for\ each\ A \in H(X)$.

Definition 2.5: A subspace M of a Banach space X is said to be (linearly) complemented (continuous linear retract) if there exists a continuous linear projection $P: X \to X$ such that P(X) = M.

(*) Equivalently, the quotient map $\varphi: X \to X/M$ admits a continuous linear lifting.

However, in the Lipschitz category only one sided implication holds. Indeed, the existence of Lipschitz right inverse for the quotient map $\varphi: X \to X/M$ yields M as a Lipschitz retract, with the retraction r being given by given by: $r(x) = x - \psi\varphi(x) + \psi(0)$. Here $\psi$ is a Lipschitz right inverse of $\varphi$. The assertion that the Lipschitz analogue of the reverse implication does not hold follows from the fact that $c_0$ is a 2-Lipschitz retract of $\ell_\infty$ (but not complemented) and that the quotient map $\varphi: \ell_\infty \to \ell_\infty/c_0$ does not admit a Lipschitz inverse (See Theorem 3.22).

2.6: Further Examples
(i) $c_0$ is a 2-Lipschitz retract of $\ell_\infty$ and that 2 is the best constant. In fact, it can be shown that any separable subspace of $\ell_\infty$ containing $c_0$ cannot be an r-Lipschitz retract of $\ell_\infty$ for $1 \leq r < 2$.

(ii) The space $C_u(M)$ of all bounded uniformly continuous functions on a metric space M (with the sup-norm) is a 12-absolute Lipschitz retract. (See [4], Theorem 1.6). For compact metric spaces K, the constant 12 has been improved to 2 by Kalton. See [25], Theorem 3.5.

(iii) For a given Banach space X, the space UC(X) of sequences in X such that the series $\sum_{n=1}^{\infty} x_n$ is unconditionally convergent is a Lipschitz retract of WUC(X), the space of sequences which are weakly unconditionally Cauchy: $\sum_{n=1}^{\infty} |f(x_n)| < \infty$, for each $f \in X^*$([27], Theorem 5.1).

(iv) $H(X)$ is an 8-abolute Lipschitzretract ([4], Theorem 1.7).

Michael's theorem motivates the following question:

**Question 2.7**: Given a Banach space X, does there exist a Lipschitz (uniformly continuous) selector $\psi: H(X) \to X$ ($\psi(A) \in A$)?

**Answer**: Yes, if X is finite dimensional.

It was conjectured for a long time that the previous question had an affirmative answer for all Banach spaces. However, the conjecture was settled in the negative by Positselskii(1971) and independently by Przeslawskii and Yost (1989).

Theorem 2.8([37], [38]): Given an infinite dimensional Banach space X, there does not exist a Lipschitz (uniformly continuous) selector $\psi: H(X) \to X$, ($\psi(A) \in A$).

Remark 2.9: It can also be deduced that for an infinite dimensional Banach space X, there exists no Lipschitz/uniform selector: $K(X) \to X$ where $K(X)$ is the sub-collection of H(X) consisting of compact sets.

The following questions arise naturally in this circle of ideas:
a. Does there exist a Lipschitz retraction: $H(X) \to X$?
b. Does there exist a uniform retraction: $H(X) \to X$?
c. Does there exist a Lipschitz retraction?
d. Does there exist a uniform retraction: $K(X) \to X$?

We have the following implications:
a. Holds if and only if X is an absolute Lipschitz retract (ALR).
b. Holds if and only if X is an absolute uniform retract (AUR).
c. Open
d. $\Rightarrow X^{**}$ is injective([38]). Converse open.

3. Lipschitz retracts in Banach space

As in topology where the theory of retracts is an important object of study, the linear counterpart of this study in the setting of Banach spaces has led to new insights into the theory and structure of Banach spaces. It turns out that the investigation of subsets of Banach spaces which arise as Lipschitz or uniform retracts of the ambient space has far reaching implications in the study of extension of Lipschitz mappings in Banach spaces. A comprehensive treatment of these ideas has been covered in [6]. A parallel development involving the extension of polynomial maps between Banach spaces has been an active area of research in recent years. See [12] for a smorsgasbord of some deep and interesting developments surrounding the extendability of polynomial maps and their implications on the geometry of Banach spaces.

We begin with the example of an absolute retract as encountered in

Proposition 2.3: A closed convex subset A of a Banach space is an absolute retract. Indeed, let M be a metric space containing A as a subset and consider the set-valued map $\phi: M \to H(X)$ where

$$\phi(x) = \begin{cases} \{x\}, x \in A \\ A, x \notin A \end{cases}$$

Clearly, $\phi$ is l.s.c. map taking values in $H(X)$, so Michael's theorem yields a continuous selection of $\phi$ which is the desired (continuous) retraction.

In Hilbert space, the situation is far more satisfactory.

Proposition 3.1: In a Hilbert space H, every closed convex subset is a nonexpansive retract (NR) of H.

(Given a closed, convex subset C of H, consider the 'metric projection':

$$P_C(x) = \{y \in C; \|x - y\| = inf(\|x - z\|; z \in C)\}.$$

The following theorem of S. Reich shows that the above property of metric projections characterises Hilbert spaces.

Theorem 3.2([39]): X is (isometrically) a Hilbert space if (and only if) each closed, convex subset of X is an (NR).

In the light of Theorem 2.2 and Reich's theorem quoted above, the following problems arise naturally:

Problem 3.3: (Lipschitz variant) Describe the class of Banach spaces X such that
(a) Unit ball of X is a Lipschitz retract of X.
(b) The above theorem (of Reich) remains true for Lipschitz retracts in place of non-expansive retracts.

Regarding (a), it turns out that the property holds in *all* Banach spaces:

Consider the (radial) mapping:

$$r(x) = \begin{cases} x, \|x\| \leq 1 \\ \dfrac{x}{\|x\|}, \|x\| \geq 1. \end{cases}$$

It follows that r maps X onto $B_X$ which is 2-Lipschitz:

$\|r(x) - r(y)\| \leq 2\|x - y\|, x, y \in X.$/Centre

However for X, a Hilbert space, we have: $\|r(x) - r(y)\| \leq \|x - y\|, x, y \in X$, i.e., r is nonexpansive. Conversely, we have:

Theorem 3.4 ([28]): X is a Hilbert space if and only if r is nonexpansive.

Remark 3.5: Regarding the existence of the metric projection in general Banach spaces, it turns out that $P_C(x)$ may not exist or if it exists, it may be multi-valued. However, in uniformly convex Banach spaces X, $P_C(x)$ exists and is unique. In the case of C being a closed convex and bounded subset of such X, it can be shown that $P_C$ is uniformly continuous on bounded neighbourhoods of C. Combining this fact with the previous considerations, we have

Theorem 3.6([24]): In a super-reflexive Banach space X, each closed, convex and bounded subset C of X is a uniform retract of X.

Indeed, assuming that C is contained in $B_X$ and that X is (in a suitable renorming of X) uniformly convex, composing $P_C$ with the 2-Lipschitz retraction r as described above, gives the desired conclusion.

The question regarding the Lipschitz analogue of Goebel's theorem is answered in the following surprising theorem of Benyamini and Sternfeld:

Theorem 3.7 ([5]): In an infinite dimensional Banach space, the unit sphere is always a Lipschitz retract of the unit ball. Equivalently, there exists a Lipchitz mapping on the unit ball of each infinite dimensional Banach space lacking fixed points.

The previous theorem suggests that in order to ensure the possible existence of a retraction with better regularity properties acting from X onto its convex subsets, it is natural to consider a subclass of the family of all closed convex subsets and wonder if the sets in this subclass may arise as retracts of mappings which may be chosen to be, say uniformly continuous or even Lipschitz continuous. Remarkable progress along these lines was achieved by Kalton (see [24]) and Hajek [16], [17] which considerably strengthens Proposition 2.3.

Theorem 3.8([16], [17], [24]): Every compact convex subset of a Banach space X is a (absolute) uniform retract of X.

The following theorem of Cheng et al shows that for a wider class of convex subsets of a Banach space, the conclusion of Theorem 3.8 is still valid.

Theorem 3.9([9]): Every convex super-weakly compact subset of a Banach space is an absolute uniform retract.

Here, *super-weakly compact* subset is meant in the following sense:

There exists $n > 1$ such that for each set of vectors $(x_i)_{i=1}^n \subseteq X$, there exists $1 \leq k \leq n$ such that $d\big(co(x_1, \ldots, x_k), co(x_{k+1}, \ldots, x_n)\big) = 0$, where $co(A)$ denotes the convex hull of $A \subseteq X$.

At this point, it may be noted that the somewhat technical definition of a super-weakly compact subset of a Banach space as given above makes sense when we note that the idea was chiefly motivated by the desire to provide a localised version of super reflexivity. In fact, analogous to the well-known fact that a Banach space is reflexive if and only if its closed unit ball is weakly compact, we have the following analogous characterisation of super reflexivity in Banach spaces.

Theorem 3.10([9]): A Banach space X is super-reflexive if and only if $B_X$ is super weakly compact.

In contrast to uniform retracts in Banach spaces as indicated above, the following theorem shows that the Lipschitz retractibility of compact convex subsets places rather severe restrictions on the underlying Banach space.

Theorem 3.11 ([17]): Every compact convex subset of X is a Lipschitz retract of X if and only if X is a Hilbert space.

As already observed in the affirmative solution to Problem 3.3(a) above which tells us that the unit ball of each Banach space X is a 2-Lipschitz retract of X, Godefroy and Ozawa [14] were led to ask the following question:

Question 3.12: Is it true that every separable Banach space X admits a compact convex generating subset M of X which is a Lipschitz retract of X? In symbols, we ask whether every separable Banach space has GCCR?

Though the question remains open in its full generality, in a series of recent papers ([16], [18], [35]), Hajek and his collaborators have shown that the following chain of implications holds:

$$Schauder\ basis \Rightarrow FDD \Rightarrow \pi-property \Rightarrow BAP \Rightarrow AP \Rightarrow CAP.$$

Here, the $\pi - property$ of a Banach space X refers to the existence of a uniformly bounded net of finite rank projections on X converging strongly to the identity map on X.

In the presence of GCCR involving compact convex sets which are 'small' in a certain technical sense - denoted GCCR* - we have the following string of implications([16], [17]):

$$FDD \Rightarrow GCCR^* \Rightarrow \pi - property \Rightarrow BAP.$$

Further, in an important work [17] the authors show that in the special case when the set involving $GCCR$ in X is 'small' in a certain technical sense, then X has the BAP.

These considerations give an indication that the answer to Question 3.12 may be negative, at least in Banach spaces lacking the bounded approximation property (BAP). However, in a recent breakthrough on a positive solution to the problem, Medina showsthat the Holder analogue of the problem admits an affirmative solution.

Theorem 3.13 ([35]): In every separable Banach space X, there exists a compact convex generating subset of X that is an $\alpha$ − Holder retract of X for each given $0 < \alpha < 1$.

In the same paper [35], the authors show that the following assertion also holds, thus providing further evidence to the previous statement that the Godefroy-Ozawa problem may have a negative solution.

Theorem 3.14: In every separable Banach space X lacking BAP and for each $0 < \alpha < 1$, there exists a compact convex generating subset of X that is an $\alpha$ − Holder retract but not a Lipschitz retract of X .

Note: At the time of resubmitting the paper after revision, an important breakthrough indicating a final push towards a negative solution to the Godefroy-Ozawa problem has been announced in a recent work by R. Medina [36]. In the said work, the author is able to construct a separable Banach space X and some $\lambda > 1$ such that no compact convex generating subset M of X admits a $\lambda$ −Lipschitz retraction onto M.

Going back to Problem 3.3(b), the next result by this author provides a complete answer which is the isomorphic analogue of Reich's theorem involving Lipschitz retracts.

Theorem 3.15 ([42]): A Banach space X has the property that each closed and convex subset of X is a Lipschitz retract of X if and only if X is isomorphic to a Hilbert space.

The above statement follows from the next theorem where for a given metric space M the symbol $Lip_0(M)$ shall be used for the space of all *Lipschitz* functions $f: M \to \mathbb{R}$, vanishing at the origin, say $\theta \in M$. It is easily checked that under pointwise operations, $Lip_0(M)$ is a Banach space when equipped with the norm:

$$\|f\|_{Lip} = \sup_{x \neq y} \frac{\sup |f(x) - f(y)|}{d(x,y)}.$$

A (closed) subset A of M which is also a Lipschitz retract of M gives rise to an extension operator $L: Lip_0(A) \to Lip_0(M)$, i.e., $L(g)|_A = g$ for each $g \in Lip_0(A)$. Bounded extension operators can also be constructed using the more general notion of a K-*random projection* which subsumes the notion of a K-Lipschitz retraction as has been done by the authors in [2]. In particular, there exists a bounded linear extension operator $Lip_0(A) \to Lip_0(M)$ for each closed subspace A of a doubling metric space M.

For later use, we shall also need to recall the notion of a *free space* $\mathfrak{J}(M)$ over a metric space M which is defined as the (unique) predual of $Lip_0(M)$: $\mathfrak{J}(M)^* = Lip_0(M)$. Concretely, $\mathfrak{J}(M)$ may be realised as the closed linear span of the set of Dirac $\delta$-functions inside $Lip_0(M)^*$:

$$\mathfrak{J}(M) = \overline{span\ \{\delta_x \in Lip_0(M)^*; \delta_x(f) = f(x), x \in M\}}.$$

Also, a subspace Y of a Banach space X is said to be *locally complemented* if there exists $c > 0$ such that for each finite dimensional subspace M of X, there exists a continuous linear map $f: M \to Y$ with $\|f\| \leq c$ and $f(x) = x$ for all $x \in M \cap Y$.

Remark 3.16: It can be shown that Y is locally complemented in X if and only if for each (closed) subspace Y of X, there exists a bounded linear 'extension' map $\psi: Y^* \to X^*$, $\psi(g)|_Y = g$ for each $g \in Y^*$ ([22], see also [19], Proposition 5.5).

Theorem 3.17: Let X be a Banach space and Y a subspace of X such that there exists a bounded linear 'extension' map $L: Lip_0(Y) \to Lip_0(X)$, i.e., $L(g)|_Y = g$ for each $g \in Lip_0(Y)$. Then L restricted to $Y^*$ gives rise to a bounded linear 'extension' map $G: Y^* \to X^*$; $G(g)|_Y = g$ for each $g \in Y^*$. Equivalently, Y is locally complemented. In particular, a Lipschitz retract is always locally complemented.

Proof (Brief sketch): The idea of the proof of the above theorem as given in [42] (see also [41]) is based on the existence of Banach limits on an abelian group which is used to produce a new bounded linear 'extension' map $G: Lip_0(Y) \to Lip_0(X)$ from L which takes values in $X^*$. In other words, the map $\psi$ restricted to $Y^*$ yields a bounded linear extension map $L: Y^* \to X^*$. Recall that a Banach limit- hereinafter to be denoted by $\int . dx$ – is a continuous linear functional on $\ell_\infty(X)$ satisfying the following conditions:

(a) $\|\int .. dx\| = 1$.
(b) $\int 1. dx = 1$.
(c) $\int f(x + x')dx = \int f(x)dx, \forall f \in \ell_\infty(X)$ and $\forall x' \in X$.

The map $G$ is now defined for $f \in Lip_0(Y)$ and $z \in X$ by the formula

$$G(f)(z) = \int \{\int [(Lf)(x+y+z) - (Lf)(x+y)]dy\}dx.$$

The linearity and boundedness of L together with (a) and the Lipschitz property of f yield that $G: Lip_0(Y) \to Lip_0(X)$ is a well-defined bounded linear map. Using translation invariance of the Banach limit as indicated in (c) and an appropriate manipulation of the Banach limit, it follows that for $f \in Lip_0(Y)$ and $z_1, z_2 \in X$ we have,

$$G(f)(z_1 + z_2) = G(f)(z_1) + G(f)(z_2).$$

and that $G(g)(z) = g(z)$ for $g \in Y^*$ and $z \in Y$. In other words, $G: Y^* \to X^*$ is an 'extension' operator. Together with Fakhoury's theorem [11], this completes the proof.

*Consequences*:
The above theorem may be looked upon as a Lipschitz analogue of the celebrated complemented subspaces theorem of Lindenstrauss and Tzafiriri characterising Hilbert spaces isomorphically in terms of complementedness of its closed subspaces.

Remark 3.18: An inspection of the proof of Theorem 3.17 shows that the adjoint map $G^*: X^{**} \to Y^{**}$ is such that for $y \in Y \subseteq X \subseteq X^{**}, G^*(y) = y$ where y is identified with the evaluation map $y^\wedge \in X^{**}$. Indeed, for $y^* \in Y^*$, we have: $\langle y^*, G^*(y^\wedge)\rangle = \langle G(y^*), y^\wedge\rangle = \langle y, G(y^*)\rangle = \langle y, y^*\rangle = \langle y^*, y^\wedge\rangle$. Thus, $G^*$ is the identity map on Y. This yields an alternative proof of the following well known theorem of Lindenstrauss [33].

Corollary 3.19: If Y is a subspace of X which is a c-Lipschitz retract and is complemented in its bidual, then Y is c-complemented in X.

Proof: Let $P: Y^{**} \to Y$ be a bounded linear projection and let $G: Lip_0(Y) \to Lip_0(X)$. be an 'extension' operator induced by a Lipschitz retraction of X onto Y. By Theorem 3.17, G may be considered as an extension operator $G: Y^* \to X^*$. Then $Q = P \circ G^*|_X$ defines a projection $Q: X \to Y$: indeed, since $Q(x) \in Y$, the previous remark gives that

$$Q^2(x) = (P \circ G^*|_X)(Q(x))$$
$$= P(Q(x)) = P^2(G^*|_X)(x) = P \circ G^*|_X(x) = Q(x).$$

Example 3.20: An application of the previous corollary to $Y = P^2(X)$, the space of 2-homogeneous polynomials on the Banach space X yields, since Y is a dual space (as the dual of a symmetric tensor product), that Z is a complemented subspace of $Lip_0(X)$ as soon as Y is a Lipschitz retract of $Lip_0(X)$. However, in certain special cases that include Hilbert spaces, it was shown by Hajek and Russo [21] that for a large class of Banach spaces that include the Hilbert space case, $P^2(X)$ is not complemented in – and is, therefore, not a Lipschitz retract of - $Lip_0(X)$. This is in sharp contrast with an old result of Lindenstrauss asserting that the space $P^1(X)(= X^*)$ of 1-

homogeneous polynomials is always 1-complemented in $Lip_0(X)$ for every Banach space X.

Remark 3.21: A Banach space X is said to be $\lambda-injective$ ($\lambda \geq 1$) if every X-valued continuous linear mapping f on any subspace of a given Banach space can be extended to a continuous linear map g on the whole space such that $\|g\| \leq \lambda\|f\|$. It can be shown that a $\lambda-$ injective Banach space is a $\lambda$- absolute Lipschitz retract ($\lambda$-ALR) i.e., it is a $\lambda$-Lipschitz retract of every metric space containing it and that an (ALR) has the property that it's an $L_1$-predual, i.e., its bidualis $\lambda-$ injective. However, an $L_1$-predual may not be an ALR. A proof of this fact depends upon the following important theorem of Kalton.

Theorem 3.22([26]): There does not exist a Lipschitz right inverse for the quotient map $\varphi: \ell_\infty \to \ell_\infty/c_0$.

Interestingly, the existence of a Lipschitz right inverse for quotient maps leads to the following useful characterisation of Hilber space among infinite dimensional Banach spaces.

Theorem 3.23([42]): Given a Banach space X such that for each closed subspace M of X and the quotient map $\varphi: X \to X/M$, there exists a Lipschitz map $\psi: X/M \to X$ with $\varphi\psi = $ identity map on $X/M$, then X is (isomorphically) a Hilbert space.

This follows because under the conditions of the theorem, the map: $r(x) = x - \psi\varphi(x) + \psi(0)$ defines a Lipschitz retraction onto M, and so weare in the situation of Theorem 3.15 to arrive at the desired conclusion.

A local theoretic analogue of the above theorem also holds as stated below.

Theorem 3.24([42]): Assume that there exists $c > 0$ such that for each finite dimensional subspace M of X the quotient map $\varphi: X \to X/M$, admits a Lipschitz right inverse Lipschitz map $\psi: X/M \to X$: $\varphi\psi = $ identity map on $X/M$ with $\|\|\psi\|\|_{Lip} \leq c$. Then X is (isomorphically) a Hilbert space.

Indeed, as was noted in the argument following Definition 2.5, the stated condition gives that each finite dimensional subspace of X is a c-Lipschitz retract. Combining this observation with Corollary 3.18 yields that finite dimensional subspaces of X are c-complemented. An application of the local version of Lindenstrauss-Tzafiriri theorem (see [1], Theorem 12.4.4) to the effect that a Banach space X is (isomorphically) a Hilbert space if (and only if) its finite dimensional subspaces are uniformly complemented completes the argument.

In view of its relevance to the theme of an earlier work by the author [41] devoted to certain aspects of nonlinear analysis in Banach spaces, Theorem

3.17 (and some of its corollaries) also appear in the said volume with proofs of some of these results having since yielded to slight modifications or improvements as included in the present paper.

Our main theorem (Theorem 3.17) motivates the following question:

Problem 3.25: Does a locally complemented subspace Y of a Banach space X imply that Y is a Lipschitz retract of X.

It turns out (see [23], Proposition 3.2) that Problem 3.25 is equivalent to one of the most important open problems in Banach space theory asking if every Banach space is a Lipschitz retract of its bidual.

*4. Extension maps involving subsets of a Banach space*

It's natural to ask if an appropriate analogue of Theorem 3.15 holds in case Y is chosen to be an arbitrary subset of the metric space (M, d). In other words, the question involves the existence of a bounded linear 'extension' operator $\psi: Lip(A) \to Lip(M)$ where A is now chosen to be an arbitrary subset of the metric space M. If such is the case, we say that the metric space M has the *'simultaneous Lipschitz extension property'* (SLEP). An answer to this question turns out to be highly nontrivial, especially in the case of X being a (non-Banach) metric space. We note that the existence of an 'extension' operator entails (i) extendability of Lipschitz mappings from arbitrary subsets to the entire space (ii) a suitable choice of the extended map which may be effected in a bounded linear manner. Whereas the validity of (i) is guaranteed by an on old result of McShane (see [40]), in the case of Banach spaces, the existence of a suitable choice of the extended map as required in (ii) forces the space to be finite dimensional as we shall see shortly. Further, the extension procedure involving (i) with the mappings consisting of contractions and taking values inside $\ell_2$ results in X being a Hilbert space ([40].Theorem 5.8).This is a well-known converse to Kirszbraun's theorem [29].

Beyond the class of Banach spaces, this question has been a subject of intensive research spanning different areas of mathematics including, in particular, geometry and group theory. There are interesting examples of situations where such extensions always exist which include, for example, $\mathbb{R}^n$ (with respect to any norm), the Heisenberg group, doubling metric spaces and metric trees (See [6]). On the other hand, making use of some deep but well known facts from the local theory of Banach spaces, it is possible to show that the extension procedure for arbitrary subsets breaks down for all infinite dimensional Banach spaces.

**Theorem 4.1**: A Banach space X has (SLEP) if and only if it is finite dimensional.

Full details of the proof may be looked up in [42].

Remarks 4.2: In contrast to Theorems 3.15 and 4.1 involving the existence/nonexistence of extension operators from spaces of Lipschitz functions on (closed) subspaces and arbitrary subsets of Banach spaces respectively, Kopecka [30] shows that in the latter case involving arbitrary subsets S of a Hilbert space H, there exist 'extension' operators $\psi: Lip_b(S,H) \to Lip_b(H,H)$ preserving the Lipschitz norm and acting between spaces of vector-valued *bounded* Lipschitz functions such that $\psi$ is continuous but *not* linear. Here the spaces of bounded Lipschitz functions are equipped with the sup-norm.

(ii) In the light of Theorem 3.15 and the corollaries following from it, the question arises regarding the extent of the validity of these conclusions under the assumption that for a subspace Z of X, the extension operator $F: Lip_0(Z) \to Lip_0(X)$ is assumed to be Lipschitz rather than linear (and continuous). Accordingly, we pose this as an open problem:

Question 4.3: Let X be a Banach space. Does the existence of a Lipschitz extension operator $F: Lip_0(Z) \to Lip_0(X)$ for each subspace Z of X yield that X is a Hilbert space?

As shown in [42], the answer is in the affirmative if X is assumed to be separable. This follows from an important theorem of Godefroy and Kalton [13] asserting that given Banach spaces X and Y with X separable and a bounded linear map $T: Y \to X$, then the existence of a Lipschitz right inverse $S: X \to Y$ of T yields the existence of a bounded linear right inverse $L: X \to Y$ of T: $ToL = Id_X$.

We conclude with the following problems which remain open.

PROBLEM 4.4: Describe sufficient conditions on a closed, convex subset C of a (separable) Banach space X such that C is a Lipschitz retract of X.

We note that by combining the fact that $\lambda(X) = \infty$ for each infinite dimensional Banach space X (Theorem 4.1) with the radial mapping $r: X \to B_X$ being a 2-Lipschitz retraction, it follows that $\lambda(B_X) = 2\lambda(X) = \infty$.

As noted in the discussion of proof of Theorem 4.1(b), an important ingredient of the proof was the following equality of Brudnyi:

$$\lambda(M) = \sup\{\nu(M,Z);\ Z \in FD\} \ldots (*)$$

In a recent work, Basso [3] shows that FD in the above equality can be replaced by the class of all dual Banach spaces. This is achieved by showing the following inequality:

$$\lambda(M,N) \leq \nu(M,N,Z), \text{ for every dual Banach space} \ldots (**)$$

This motivates the following question.

PROBLEM 4.5: Is it true that the equality (*) holds for the class of all Banach spaces? In particular, is it possible to replace FD by the class of all Banach spaces which are complemented in their bidual?

In a related work that is in progress, it has been observed that a closer look at the proof of Basso's theorem actually yields an affirmative answer to a strengthened form of the previous problem involving Banach spaces which are 1-Lipschitz retracts of their biduals.

In the same work [3], Basso shows how $\lambda(M)$ is related to ae(M), the absolute Lipschitz constant $ae(M)$ of the metric space M:

$$ae(M) = \sup\{v(M,N); M \subset N\}.$$

Recalling the notation

$$\lambda(M,N) = \inf\{\|T\|; T \in \text{Ext}(M,N)\}$$

we see that (**) yields:

Theorem 4.6([3]): For a metric space M, we have

$$ae(M) \geq \sup_{N \supset M} \lambda(M,N),$$

where the equality holds for finite metric spaces M and finite metric spaces N (containing M).

Corollary 4.7(a) ([3]): $ae(3) = 3/4$.
(b) $\dfrac{5+4\sqrt{2}}{7} \leq ae(4) \leq \dfrac{3+6\sqrt{2}}{7}$.

We conclude with a brief mention of a recent breakthrough on a nonlinear version of Argyros's discovery of (hereditarily) indecom-posable Banach spaces: there exist Banach spaces X which do not admit (proper) infinite dimensional complemented subspaces. In other words, X does not admit (linear) continuous projections onto closed infinite dimensional subspaces of X.

The following theorem of Hajek and Quilis provides a nonlinear analogue of the previous statement.

Theorem 4.8 ([20]): There exists a complete metric space M with $|M| = c$ such that no non-singleton closed separable subspace of M is a Lipschitz retract of M.

*Some more examples of retracts in analysis*

In the linear theory of Banach spaces, it's well known, and easy to prove, that $K(\ell_2)$, the space of compact operators is not complemented in

$L(\ell_2)$. It's natural to ask if $K(\ell_2)$ may be a Lipschitz retract of $L(\ell_2)$. In a recent paper, R. Tanaka[44] shows that, more generally, every 2-sided ideal in a Von Neumann algebra is a Lipschitz retract in the ambient space.

In the more general case of $p \geq 1$ it's possible to show that $K(\ell_p)$ is not complemented in $L(\ell_p)$ whereas the two classes of operators $K(\ell_p, \ell_q)$ and $L(\ell_p, \ell_q)$ coincide (by Pitt's theorem) if $\infty > p > q \geq 1$. Surprisingly, regardless of the range of p and q, the Lipschitz analogue of this phenomenon always holds, as the following theorem shows.

Theorem 4.9([9]): For $1 \leq p, q < \infty$, the space $K(\ell_p, \ell_q)$ of compact operators is a 8-Lipschitz retract of $L(\ell_p, \ell_q)$.

The previous result may be seen as a far reaching strengthening of Example 2.6(ii).

*Symmetric products* (Kovalev theory)

By definition, an injective Banach space is complemented in each Banach space containing it as a subspace. The metric version of this property is provided by the notion of hyperconvexity which is characterised by the property that a hyperconvex metric space is an absolute Lipschitz retract (ALR) - a Lipschitz retract of each larger metric space. The Hahn Banach theorem yields that $\mathbb{R}$ is injective as a Banach space. As a metric space, $\mathbb{R}$ is also hyperconvex (ALR).

In an important work, Kovalev [31] shows that the symmetric product $\mathbb{R}^{(n)}$ of $\mathbb{R}$ is also an ALR for all $n \geq 1$. Here, $\mathbb{R}^{(n)}$ is defined to be the set of nonempty subsets of $\mathbb{R}$ with cardinality at most n equipped with the Hausdorff metric. The same also holds for higher dimensional Euclidean space $\mathbb{R}^d$. This follows from his main theorem in a subsequent work [32], asserting that for all $1 \leq k \leq n$ there exists a Lipschitz retraction from $H^{(n)}$ onto $H^{(k)}$ where H is a Hilbert space, finite or infinite dimensional. Here the proof is based on the existence of gradient flow trajectories in a finite dimensional subspace which is guaranteed by classical theory of ordinary differential equations. Another delightful consequence of the latter assertion is that $H^{(n)}$ embeds in a biLipschitz manner into an Euclidean space. In fact, a metric space X admits a bi-Lipschitz embedding into an Euclidean space if and only if its symmetric product $X^{(n)}$ does. From this, it can be deduced that $\mathbb{R}^{(n)}$ is an ALR. However, the author also remarks that it remains unknown if the property of being an ALR passes from X to its symmetric product $X^{(n)}$. The question whether for a Banach space X, $X^{(k)}$ is a Lipschitz retract of $X^{(n)} (1 \leq k < n)$ remains wide open.

M. A. Sofi  
Department of Mathematics  
J. K. Institute of Mathematics  
(Affiliated to Kashmir University)  
Srinagar-190008  
India  
email: aminsofi@gmail.com